\documentclass[a4paper, 12pt]{amsart}
\usepackage{amsgen,amsmath,amstext,amsbsy,amsopn,amsfonts,amssymb}
\newtheorem{theorem}{Theorem}
\newtheorem{lemma}[theorem]{Lemma}
\newtheorem{corollary}[theorem]{Corollary}
\newtheorem{proposition}[theorem]{Proposition}
\newtheorem{obs}[theorem]{Observation}

 \newtheorem{defi}[theorem]{Definition}
\newenvironment{definition}{\begin{defi}\rm}{\end{defi}}
\newtheorem{exa}[theorem]{Example}
\newenvironment{example}{\begin{exa}\rm}{\end{exa}}
\newtheorem{rem}[theorem]{Remark}
\newenvironment{remark}{\begin{rem}\rm}{\end{rem}}
\newtheorem{rems}[theorem]{Remarks}

\newtheorem{ack}[theorem]{Acknowlegment}

\def\proof{\noindent\textbf{Proof}\quad}
\def\bsq{\blacksquare\medskip}
\def\n{\noindent}
\def\H{\mathcal H}
\def\L{\mathcal L}
\def\K{\mathcal K}
\def\P{\mathcal P}

\def\C{\mathcal C}

\def\R{\mathcal R}
\def\L{\mathcal L}

\def\CalS{{\rm Sym}}

\def\NN{{\mathbf N}}
\def\ZZ{{\mathbf Z}}
\def\CCC{{\mathbf C}}
\def\RRR{{\mathbf R}}
\def\QQ{\mathbf Q}
\def\RR+{{\mathbf R}^*}

\def\Q_p{{\mathbf Q}_p}
\def\S1{{\mathbf S}^1}
\def\OO{{\mathbf O}}

\def\Ind{{\rm Ind}}

\def\PT{Property $(T)$}
\def\PF{Property $(F)$}

\def \Tlp {Property $({T}_{\ell_p})$}
\def\PTB{Property $({T}_B)$}
\def\Ptau{Property $(\tau)$}
\def\eps{\varepsilon}
\def\Ga{\Gamma}
\def\ga{\gamma}
\def\la{\lambda}
\def\vfi{\varphi}

\def\tous{\qquad\text{for all}\quad}

\begin{document}

\title{On groups with  Property $({T}_{\ell_p})$}
\author{Bachir Bekka}
\address{Bachir Bekka \\ IRMAR \\ UMR-CNRS 6625 Universit\'e de  Rennes 1\\
Campus Beaulieu\\ F-35042  Rennes Cedex\\
 France}
 
\email{bachir.bekka@univ-rennes1.fr}
\author{Baptiste Olivier}
\address{Baptiste Olivier \\ IRMAR \\ UMR-CNRS 6625\\ Universit\'e de  Rennes 1\\
Campus Beaulieu\\ F-35042  Rennes Cedex\\
 France}
\email{baptiste.olivier@univ-rennes1.fr}

\thanks{The first author acknowledges the support of the  French Agence Nationale de la Recherche (ANR) for the projects ANR-2010-BLANC
 GGAA and  ANR-2009-BLANC AGORA}

\begin{abstract}
Let $1<p<\infty, p\neq 2.$  \Tlp\ for a 
second countable locally compact  group $G$
  is a weak  version of Kazhdan \PT, defined in terms of 
  the orthogonal representations of $G$ on $\ell_p.$
    \Tlp \ is characterized by 
an isolation property of the trivial representation
from the monomial unitary   representations of $G$
associated to open subgroups.
Connected groups with   \Tlp \  
are the connected groups with a compact abelianization.

In the case of a totally disconnected group,
 isolation of the trivial representation
from the quasi-regular representations 
associated to open subgroups suffices to characterize \Tlp.
Groups  with \Tlp \ 
share some important properties with Kazhdan groups
(compact generation,  compact abelianization, ...).
Simple algebraic groups over non-archimedean
local fields as well as  automorphism groups of $k$-regular trees 
for $k\geq 3$ have  \Tlp .

 In the case of discrete groups,   \Tlp\ implies  Lubotzky's \Ptau\ 
 and is implied by  \PF \ of Glasner and Monod.
%the reverse implications
%do not hold in general.
We show that an irreducible lattice $\Ga$
in a product $G_1\times G_2$ of locally compact  groups  have \Tlp,
whenever $G_1$ has \PT\
and $G_2$ is connected and minimally almost
periodic.  Such  a lattice does not have \PT\
if $G_2$ does not have \PT.
% Other examples of discrete groups with \Tlp  are 
 %groups which have the Glasner and Monod \PF.
 %Some of them do not have \PT.

\end{abstract}
\maketitle
\section{Introduction}

A fundamental rigidity property of groups, with a wide range of applications,
  is the by now classical  Kazhdan's   \PT\ defined 
in \cite{Kazhdan}.  
Bader, Furman, Monod, and Gelander  introduced in \cite{BFGM}
 \PTB\ for a general Banach space $B$, the case where 
 $B$ is a Hilbert space  corresponding to Kazhdan's  \PT.
 It was shown in  \cite[Theorem 1]{BFGM}  that $T_{L_p([0,1])}$ 
 for some $1<p<\infty$ is equivalent to Kazhdan's  \PT,
 at least  for  second countable locally compact groups.
 Moreover, \PT\ implies $(T_{L_p(X,\mu)})$  for any $\sigma$-finite measure
 $\mu$ on a standard Borel space $X$ and any $1\leq p<\infty.$

 In this article, we study \Tlp\ for $1<p<\infty$ and $p\neq 2,$ where $\ell_p$ is the usual  Banach space
 $ \ell_p (\NN)$  of $p$-summable {complex-valued}  sequences.
 (We prefer to work with \emph{complex} Banach spaces, as this is more suited to our approach
 which is based on unitary representations on complex Hilbert spaces.)
 
  Let   $G$ be a second countable  locally compact group.
  We recall the definition from \cite{BFGM} of  \PTB\ for $G$ in the  special case
  $B=\ell_p.$
  
 Denote by $\OO(\ell_p)$
  the group of linear bijective isometries of $\ell_p.$   An orthogonal representation  of $G$ on  $\ell_p$ is a homomorphism of 
   $\pi: G\to \OO(\ell_p),$ which is continuous in the sense  that
the mapping $G\to \ell_p, g\mapsto \pi(g) f$ is continuous for every $f\in \ell_p.$
A  sequence $(f_n)_n$ in $\ell_p$  is
said to be a sequence of almost  vectors if 
$\Vert f_n\Vert=1$ and
$$
\lim_n \Vert \pi(g)f_n-f_n\Vert = 0, \quad \text{uniformly on compact subsets of} \ G.
$$

Let $\pi^*: G\to \OO(\ell_q)$ be the dual representation of $\pi$ on $\ell_q,$ where $q$ is the conjugate exponent
of $p.$ Denote by  $\ell_p^{\pi(G)}$ and  $\ell_q^{\pi^*(G)}$ the closed subspaces of $G$-invariant vectors in 
$\ell_p$ and $\ell_q$, respectively.
Let $\ell'_p(\pi)$ be the annihilator of $\ell_q^{\pi(G)^*}$  in $\ell_p.$ Then we  a have a decomposition
into $G$-invariant subspaces   $\ell_p= \ell_p^{\pi(G)}\oplus \ell_p'(\pi)$.
The group $G$ has \Tlp\ if  there exists no sequence of almost  invariant 
vectors in $ \ell'_p(\pi).$ 

As we will show (see Section~\ref{S:Banach}), 
Banach's description from  \cite[Chap. XI]{Banach}
of  the group $\OO(\ell_p)$  implies 
that the orthogonal representations  of a group $G$ on  $\ell_p$
have a simple structure for  $p\neq 2$ .

In order to state our first result,  we   introduce  Property $(T; \R)$ 
with respect to a set $\R$ of unitary representations of $G;$ our definition can  easily be
reformulated in terms of the Property  $(T;\R)$  introduced in \cite[Definition 1.1]{Lubotzky-Zimmer}
for sets of irreducible unitary representations of $G.$
\begin{definition} 
\label{Def-TRep}
Let $\R$ be a set of (equivalence classes of)  unitary representations of
$G$ on Hilbert spaces. We say that $G$ has Property $(T; \R)$ if,  the trivial representation
$1_G$ of $G$ is  not  weakly contained the  direct sum $\oplus_{\pi\in \R} \pi',$
where $\pi'$ denotes the restriction of $\pi$ to the orthogonal complement  of 
the $\pi(G)$-invariant vectors in the Hilbert space of $\pi$
(that is, if the unitary representation  $\oplus_{\pi\in \R} \pi'$ has  no sequence of almost  invariant vectors).
\end{definition}
Observe that  \PT \ for $G$ corresponds to the case where $\R$ is the set of 
all (equivalence classes of)  unitary representations of $G.$
% Interesting examples of groups which are not Kazhdan groups and have $(T; \R)$ for some natural sets $\R$
%of representations are given in \cite{Lubotzky-Zimmer}.

%Let $H$ be an closed  subgroup of $G$ 
%such that the homogeneous space $G/H$ has a $G$-invariant
%Borel measure. The quasi-regular representation 
%$\la_{G/H}$ is the natural unitary representation
%of  $G$ on $L_2(G/H).$
%Recall that  $la_{G/H}$ is unitarily equivalent to the induced representation
%$\Ind_H^G 1_H$,
%where $H$ is a closed subgroup of $G$ and $\chi: H\to \S1$ is a unitary character of $H.$
%Examples  of monomial representations are   the quasi-regular representations $\la_{G/H}$ of $G$ on $L_2(G/H)$
%attached to  closed subgroups $H,$ which  correspond to the case $\chi=1_H.$ 
Recall that a unitary representation $\sigma$ of  $G$ is \emph{monomial}
if $\sigma$ is unitarily equivalent to an induced representation
$\Ind_H^G \chi,$
where $H$ is a closed subgroup of $G$ and $\chi: H\to \S1$ is a unitary character of $H.$
Examples  of monomial representations are   the \emph{quasi-regular} representations $\la_{G/H}$ of $G$ on $L_2(G/H)$
and correspond  to the case where $\chi=1_H.$ 

\begin{theorem}
\label{Theo1}
  Let   $G$ be a second countable  locally compact group. The 
  following properties are equivalent.
  
  \n
  (i) $G$  has    \Tlp \ for some $1< p<\infty, p\neq 2.$ 
%\item $G$  has    \Tlp  for every  $1<p<\infty, p\neq 2.$ 

\n
(ii) $G$ has Property $(T; \R_{\rm mon})$,
where $ \R_{\rm mon}$ is the set of monomial 
unitary representations  $\Ind_H^G \chi,$ 
associated to \emph{open} subgroups $H$ of $G.$
\end{theorem}

In particular, if $G$ has  \Tlp \  for some $1< p<\infty, p\neq 2,$
then $G$ has  \Tlp \ for every $1< p<\infty, p\neq 2.$
It is also clear from the previous theorem that \PT \ implies \Tlp; this is also a special
case of  Theorem~A in \cite{BFGM}, as mentioned above. 
As we will see later, \Tlp \ is strictly weaker than \PT.

When  $G$ is connected, the only open subgroup
of $G$ is $G$ itself and   $ \R_{\rm mon}$ 
therefore coincides with the group of unitary characters of $G,$
that is, with the Pontrjagin dual of the abelianization $G/\overline{[G,G]}.$
The following corollary is therefore an immediate consequence of Theorem~\ref{Theo1}.
 \begin{corollary}
\label{Cor0}
A connected locally compact second countable group $G$  has \Tlp \ for   $1< p<\infty, p\neq 2,$
if and only if its abelianization $G/\overline{[G,G]}$ is compact.
\end{corollary}

Our next result shows that, when $G$ is totally disconnected,  isolation of $1_G$
in the set of quasi-regular representations associated to open subgroups
suffices to  characterize \Tlp\ for $p\neq 2.$

\begin{theorem}
\label{Theo2}
  Let   $G$ be a  totally disconnected, second countable  locally compact  group. The 
  following properties are equivalent.
  
  \n
  (i) $G$  has    \Tlp \ for some $1< p<\infty, p\neq 2.$ 
%\item $G$  has    \Tlp  for every  $1<p<\infty, p\neq 2.$ 

\n  
(ii) $G$ has Property $(T; \R_{\rm quasi-reg})$,
where $ \R_{\rm quasi-reg}$ is the set of quasi-regular
representations  $(\la_{G/H},\ell_2(G/H))$ 
associated to \emph{open} subgroups $H$ of $G.$

\end{theorem}

\begin{remark}
\label{Rem-Theo2}
(i) Observe that the previous theorem does not hold, in general,
if $G$ is not totally disconnected,
as already the example $G=\RRR$ shows.

\n
(ii) The result in the previous theorem can be rephrased in terms of the existence
of an appropriate Kazhdan pair $(Q,\eps)$ as in the case of 
\PT: 
%(see   \cite[Chap. I, \S 1]{BHV}):
$G$  has    \Tlp \ for some $p\neq 2$ if only if there exists
a compact  subset $Q$ of $G$ and an $\eps>0$ such that
$$
\sup_{g\in Q}\Vert \la_{G/H}(g)f-f\Vert \geq \eps 
$$
for every open subgroup $H$ and every unit vector $f$ in the orthogonal complement
of the space of $G$-invariant vectors in $\ell_2(G/H).$

\n
(iii) We will give below (Example~\ref{Exa-Tau} and \ref{Exa-F})  examples of discrete groups  with
\Tlp\ and without \PT; these examples, which seem to be the first
of this kind,  show that isolation of trivial representation in the 
family of all quasi-regular representations does not suffice in order to
imply \PT.
\end{remark}

Groups  with \Tlp\
share some important properties with Kazhdan groups.

\begin{theorem}
\label{Theo-Consequences}
 Let $G$ be second countable
 locally compact group.   
Assume that $G$ has \Tlp\ for some $1<p<\infty$.
The following statements hold:

\n
(i) $G$ is compactly generated.

\n
(ii) The abelianized group $G/\overline{[G,G]}$ is compact.

\n
(iii)  Every  subgroup of finite index in $G$ and every 
topological group containing $G$ as a finite index subgroup 
has \Tlp. (In other words, \Tlp\ only depends on the commensurabilty class of $G.$)

\n
(iv) If amenable and totally disconnected, then $G$ is compact.
\end{theorem}

\begin{remark}
\label{Rem-Lattices}
(i) It follows from
the previous theorem that, for instance,
(abelian or non-abelian) free groups  as well as the groups $SL_n(\QQ)$ do not have 
\Tlp. 

\n
(ii) 
 \Tlp\  for $p\neq 2$
is not inherited by lattices, even in the totally disconnected case.
 Indeed,  $SL_2(\QQ_l)$ has \Tlp\  for $p\neq 2$
(see Example~\ref{Exa1} below), whereas
torsion-free discrete subgroups in $SL_2(\QQ_l)$ are free groups
(see Chap. II,  Th\'eor\`eme 5 in  \cite{Serre}).
 \end{remark}

The next  result will provide us with a class of examples of
 totally disconnected non discrete groups
with \Tlp\  for $p\neq 2$ and without \PT. 

A locally compact group $G$
has the Howe-Moore property if, for every 
unitary representation $\pi$ of $G$ without non-zero
invariant vectors, the matrix coefficients of $\pi$ are in $C_0(G).$  
For an extensive study of groups with this property, see \cite{CCLTV}. 

\begin{theorem}
\label{Theo-HoweMoore}
Let $G$ be a totally disconnected group, second countable locally compact
group with the Howe-Moore property.
Assume that $G$ is non-amenable.
Then  $G$ has  \Tlp\ for every  $1\leq p<\infty, p\neq 2.$
\end{theorem}

\begin{example}
\label{Exa1}
\n
(i)  Let $k$ be a non archimedean local field, $\mathbb G$ a simple linear algebraic 
group over $k$ and $G= {\mathbb G}(k)$  the group of $k$-points
in  $\mathbb G$ (an example is  $G=SL_n(\QQ_l)$ for $n\geq 2,$ where $\QQ_l$ is the field of $l$-adic numbers
for a prime $l$).
Then $G$  has the Howe-Moore property (see Theorem 5.1 in \cite{Howe-Moore}). 
Moreover, $G$ is amenable if and only if $G$ is compact.
So, $G$   has \Tlp \ for $p\neq 2.$
Observe that, if $k - {\rm rank} (\mathbb G) =1$,  
then $G$ does not have \PT; see Remark 1.6.3 in \cite{BHV}.
 This is, for instance, the case  for $G= SL_2(\QQ_l).$

\n
(ii) Let $G={\rm Aut}(T)$ be the group of  color preserving 
automorphisms of a $k$-regular 
tree $T$ for $k\geq 3$ or of a bi-regular tree of type $(m,n)$ for 
$m,n\geq 3.$ Then  $G$ is a totally disconnected locally
compact group and,  as  shown in \cite{Lubotzky-Mozes},
$G$ has the Howe-Moore property.
Since $G$ is non-amenable, it has  \Tlp\ for $p\neq 2$. Observe that  $G$ 
does not have \PT. 
\end{example}

We turn now  to  discrete groups.  
The examples  of discrete groups with  \Tlp\ for $p\neq 2$ and without \PT\
we give below 
are   related either to Lubotzky's \Ptau\  or  to  the Glasner-Monod \PF.

Recall that a discrete group $\Ga$
has \Ptau, if $\Ga$ has Property $(T; \R_{fi})$,
with respect to the set  $\R_{fi}$ of regular representations $\la_{\Ga/H}$
associated to subgroups $H$ of finite index (\cite[Definition 4.3.1]{Lubotzky}). The following result
is an immediate consequence of  Theorem~\ref{Theo2}.

\begin{proposition}
\label{Pro-Tau}
Every discrete group with \Tlp\
has  Lubotzky's  \Ptau. $\bsq$
\end{proposition}

\begin{remark}
\label{Rem-Tau0}
\Ptau \  does not imply \Tlp\ for $p\neq 2.$
Indeed, by Theorem~\ref{Theo-Consequences}, every countable discrete group  with 
\Tlp\ is finitely generated. However, there are 
 groups with  \Ptau\   which are not  finitely generated.
An example of such a group is  $\Ga=SL_n(\ZZ[1/\P])$ for $n\geq 2,$
where  $\P$ is an infinite  set of primes not containing all primes and 
$\ZZ[1/\P]$ denotes the ring of rational numbers whose denominators
are only divisible by primes from $\P$ (see the remarks after Corollary 2.7 in  \cite{Lubotzky-Zimmer}).
\end{remark}

%In the sequel, we will give  examples of finitely generated groups 
%with  Property $(\tau)$ and without \Tlp.

We  prove \Tlp\  for the following class of  lattices.
 %in semisimple groups, which arise as follows. Let $k_1$ and $k_$ be local fields, $\mathbb G_1$ and 
%$\mathbb G_2$ be  simple linear algebraic 
%group over $k_1$ and$k_2$ and let  $G= {\mathbb G_1}(k_1)\times \mathbb G_2}(k_1).$ 
%Let $\Ga$ be an irreducible lattice in $G.$ Assume that ${\mathbb G_1}(k_1)$ has \PT (this is
%exactly  the case when $k_1$-{\rm rank}\( \mathbb G_1) \geq 2$, see \cite{BHV}).
Recall that a lattice $\Ga$ in a product $G_1\times G_2$ of 
locally compact  groups is irreducible if the natural  projections 
of $\Ga$ to $G_1$ and $G_2$ are dense.
A locally compact group 
is  minimally almost periodic if it has no non-trivial 
finite dimensional  unitary representation.
\begin{theorem}
\label{Theo3}
Let $G_1, G_2$ be  locally compact second countable groups
and  $\Ga$  an irreducible lattice in  $G=G_1\times G_2.$
Assume that $G_1$ has \PT\ 
 and that 
 $G_2$ is connected and minimally almost periodic.
Then $\Ga$ has  \Tlp\ for $p\neq 2.$
\end{theorem}
\begin{remark}
\label{Rem-Tau}
(i)   Property $(\tau)$ was established  in  \cite[Corollary 2.6]{Lubotzky-Zimmer}
 (see also Corollary of Theorem~A in  \cite{Bekka-Louvet})
 for the groups $\Ga$
appearing in Theorem~\ref{Theo3}, without the connectedness assumption 
on $G_2$.

\n
(ii) If  $G_2$ does not have \PT\, then neither does $\Ga$, since \PT\ is inherited by lattices
(\cite[Theorem~1.7.1]{BHV}).  

\end{remark}
\begin{example}
\label{Exa-Tau}
Examples of groups $\Ga$ 
as in Theorem~\ref{Theo3} are, for instance, lattices in $SO(n, 2)\times SO(n+1,1)$ for $n\geq 3.$
(Observe that Theorem~\ref{Theo3} still applies  when the connected component 
of $G_2$ has finite index, in view of  Theorem~\ref{Theo-Consequences}.iii.)
Such lattices can be obtained by the following well-known construction.
Let $\mathbf G=SO(q)$ be the orthogonal group of the quadratic
form $q(x) = x_1^2+ \cdots+ x_n^2 -x_{n+1}^2-\sqrt{2} x_{n+2}^2.$
Then ${\mathbf G}(\RRR) \cong SO(n,2)$.
Let $\sigma$ be the non trivial field automorphism of $\QQ(\sqrt{2})$ and 
${\mathbf G}^\sigma=SO(q^\sigma)$  the orthogonal group 
of  the conjugate form $q^\sigma$. Then 
 ${\mathbf G^\sigma}(\RRR) \cong SO(n+1,1)$
 and $\Ga={\mathbf G} (\ZZ[\sqrt 2])$ embeds as an irreducible lattice in ${\mathbf G^\sigma}(\RRR)\times {\mathbf G}(\RRR)$
 by means of the mapping $ \ga\mapsto (\ga^\sigma, \ga).$
 \end{example}

Another class of groups with  \Tlp\  are the groups  with \PF\ of Glasner and Monod.
We first recall how this later property  is defined.

 A continuous action of a locally compact group $G$  on a discrete countable space  $X$ is said to be amenable  
if the natural representation  of $G$ on  $\ell_2(X)$ weakly contains $1_X.$
One should mention that this notion, extensively studied in \cite{Eymard} in the case of (non necessarily discrete) homogeneous spaces,
 is different from Zimmer's notion of amenable group action from Section~4.3 in \cite{Zimmer}.
A  discrete group $\Ga$ has  \PF\ if
every  amenable action of  $\Ga$ on a  countable space  $X$
has a fixed point (\cite[Definition 1.3]{Glasner-Monod}).
It turns out that \PF\  implies \Tlp \  for  $p\neq 2.$

\begin{proposition}
\label{Pro-F}
Let  $\Ga$ be a discrete countable group with the Glasner-Monod \PF. Then   
$\Ga$ has  \Tlp\   for  $p\neq 2.$
\end{proposition}

\begin{remark}
\label{Rem-F-Tlp}

 \Tlp\  does not imply \PF. 
Indeed,  a group with \PF \  has no non-trivial finite quotient. 
However there are groups, such as  $SL_3(\ZZ)$,
which have   \PT \, and hence  \Tlp,  and which have non-trivial finite quotients.

%It would be interesting to know whether \Tlp\  for $p\neq 2$ and \PFO\  are equivalent.
%We suspect that this is  not the case.
\end{remark}

\begin{example}
\label{Exa-F}
 It follows from Proposition~\ref{Pro-F}   that the examples  given in \cite{Glasner-Monod} of groups
with  \PF\ and without \PT\  are at the same time examples of 
groups with \Tlp\ for $p\neq 2$ and without \PT. We briefly  recall their construction.

A remarkable feature of  \PF\
is that the class of groups with this property  is preserved by finite free products
(Lemma 3.1 of \cite{Glasner-Monod}).
The free product $\Ga$ of  two non-trivial groups with \PF\ 
therefore has   \Tlp\ for $p\neq 2.$ Observe that $\Ga$  does not have \PT,
since $\Ga$ acts on a tree without fixed point.
 It remains to give examples of groups  with \PF.
 The examples given in \cite{Glasner-Monod} are  infinite 
simple Kazhdan groups or non amenable 
groups for which every proper subgroup is finite (``Tarski monsters").
 Examples of  the first kind of groups 
 were constructed by  Gromov (Corollary 5.5.E in \cite{Gromov}) as quotients of  hyperbolic groups
with \PT; for another construction, see Corollary~21 in \cite{Caprace-Remy}. 
 Examples of the second kind of groups were given by   Ol'shanskii \cite{Olshanskii}.
 \end{example}
 This paper is organized as follows.
 In Section~\ref{S:Banach} contains some basic remarks on the structure of 
 orthogonal group representations on $\ell_p$ for $p\neq 2.$
  Theorems~\ref{Theo1} and \ref{Theo2} are proved in Section~\ref{S:Theo1-2}
  and Theorems~\ref{Theo-Consequences} and \ref{Theo-HoweMoore} in Section~\ref{S:Consequences}.
 Section~\ref{S:Discrete} is devoted to  discrete groups with \Tlp\ and contains
 the proof of Theorem~\ref{Theo3} and Propositions~\ref{Pro-Tau} and \ref{Pro-F}.

 After this work was completed, we learned of the  preprint \cite{Cornulier-FM} which contains some
 overlap with our present work. Consider the following weaker version of \PF, introduced
(without name) in \cite[Remark~1.4]{Glasner-Monod}: 
a  group $G$  has Property FM, if every  amenable action of  $G$ on a  countable space 
has a finite orbit.  Using Theorem~\ref{Theo2}, it is easy to see that 
 a group has \Tlp \ if and only if it has Property FM and \Ptau. 
 Property FM is  studied  \cite{Cornulier-FM},
in an independent way and with a different motivation. It is shown there that  the lattices $\Ga$ appearing in 
our Theorem~\ref{Theo3} have Property FM. As these groups are known to have \Ptau, this gives a different proof of \Tlp\ for these lattices.
 Moreover, in connection with our Example~\ref{Exa1}.i, it  is  proved in \cite{Cornulier-FM} that semisimple algebraic groups 
 over non-archimedean local fields have Property FM.
 
   \medskip
\noindent
\textbf{Acknowlegments} We are grateful to  Yves de Cornulier, Vincent Gui\-rar\-del, Pierre de la Harpe, Nicolas Monod, and Mikael  de la Salle for useful discussions.

 \section{Orthogonal representations on $\ell_p$ for $p\neq 2$}
 \label{S:Banach}
 We begin with some preliminary remarks on 
  permutation representations  of topological groups 
  twisted by a cocycle with values in $\S1.$
  Let $G$ be a  topological group. Let  $X$ be a  
  discrete space equipped with a $G$-action. We assume that this action is continuous,
  or, equivalently, that the   stabilizers of points in $X$ are open subgroups of $G.$
  Let $c: G\times X\to \S1$ be a continuous cocycle with values in $\S1;$
  thus, $c$ 
  satisfies the cocycle relation:
$$
c(g_1g_2, x) = c(g_1, g_2x)c(g_2,x),  \tous  g_1,g_2\in G,\ x\in X. \eqno{(*)}
$$
We associate to the $G$-action and the cocycle $c$
the permutation  representation  twisted by $c,$
which is the continuous representation of $G$ on $\ell_2(X),$
denoted by $\la_{X}^c$ and
defined by the formula
$$
\la_{X}^c(g)(f)(x) = c(g^{-1},x) f(g^{-1}(x)),
\tous  g\in G, f\in\ell_2(X),  x\in X.
$$

The following lemma is a very special case of  Mackey's imprimitivity theorem
(see Theorem~3.10 in \cite{Mackey}).

\begin{lemma}
\label{Imprimitivity}
Assume that $G$ acts transitively on $X$.
Let $x_0\in X$ and denote by $H$ the stabilizer  of $x_0$ in $G.$
Let $\chi: H\to \S1$ be defined by $\chi(h) = c(h, x_0)$ for all $h\in H.$
Then  $\chi$ is a unitary character of $H$ and $\la_{X}^c$ is unitarily
equivalent to the monomial representation $\Ind_H^G\chi.$
\end{lemma}

\proof
The fact that $\chi$ is a homomorphism follows immediately from
the cocycle relation $(*).$ 

Fix a set $T\subset G$ of representatives
for the left cosets of $H.$  
The space  $\ell_2(X)$ is the direct sum $\oplus_{x\in X} V_{x}$,
where  $V_x$ is the one-dimensional space $\CCC\delta_x.$
The restriction of $\la_{X}^c$ to $H$ 
leaves  $V_{x_0}$ invariant,  with the corresponding $H$-representation
given by the character $\chi.$ Moreover,
we have  $\la_{X}^c(t) V_{x_0}= V_{tx_0}$ for all $t\in T.$
This shows that  $\la_{X}^c$ is equivalent to $\Ind_H^G\chi,$
by the defining property  of   induced representations
(see Section~3.3 in \cite{Serre-Representations};  the
case of finite groups  treated there carries over to induced representations 
from open subgroups of infinite groups).
$\bsq$

\begin{remark}
\label{Rem-MonomialCocycle}
Conversely, as is well-known, every monomial
representation of $G$ associated to an open subgroup $H$
can be realized as a  representation of the form 
$\la_{X}^c$ for the action  of $G$ on $X=G/H$
and a continuous cocycle  $c: G\times G/H\to \S1$
We recall briefly the construction of $c.$  Choose  a
 section $s: G/H\to G$ for the canonical projection $p: G\to G/H,$
 with $s(H)= e.$
 Define a
a cocycle $\alpha:G\times X\to H$ with values in $H,$ given by
$$
\alpha(g, x) = s(gx)^{-1} gs(x), \tous g\in G, x\in X.
$$
Then $c: G\times G/H\to \S1$ is defined by $c(g,x)= \chi(\alpha(g, x)).$
\end{remark}

The following lemma is an immediate consequence
of Lemma~\ref{Imprimitivity}
and the previous remark.
\begin{corollary}
\label{Cor-DirectSumMonoRep}
 Let $G$ be a  topological group.
 
\n
(i)  Let  $X$ be a   discrete space equipped with a continuous $G$ action
and  let $c: G\times X\to \S1$ be a continuous cocycle.
The associated representation $\la_{X}^c$  of $G$ on $\ell_2(X)$
is  equivalent to a direct sum of monomial representations associated to 
open subgroups of $G.$

\n
(ii) Let $\pi= \oplus_{i\in I} \Ind_{H_i}^G\chi_i$ be a  a direct sum of monomial representations associated to 
open subgroups $H_i$ of $G.$ 
Set $X=\coprod_{i\in I}G/H_i,$ the disjoint sum of the $G/H_i$'s,
with the  obvious $G$ -action. Then $\pi$ is equivalent 
to the representation $\la_{X}^c$  of $G$ on $\ell_2(X)$
for a cocycle $c: G\times X\to \S1.\bsq$ 
 \end{corollary}

 Next, we study   orthogonal  representations  of topological groups on $\ell_p$ 
  for $p\neq 2.$  
  We first recall Banach's  description of $\OO(\ell_p)$ for $1\leq p<\infty, p\neq2$
 from Chapitre~XI in \cite{Banach}.
 
 Let  $X$ be an infinite countable set and $\CalS(X)$ the
group  of all  permutations of $X.$
Let $U\in \OO(\ell_p(X))$. There exists a unique permutation 
$\sigma\in \CalS(X)$ and a unique function $ h: X\to \S1$ 
such that 
$$
U(f)(x) = h(x) f(\sigma(x)),
\tous  f\in\ell_p(X),  x\in X.
$$
(One should observe that  Banach's theorem is stated  in \cite{Banach} for 
spaces of \emph{real-valued} sequences; however,
the arguments remain valid for complex-valued sequences and yield the result stated above.)

Let  $G$ be a topological group  and $\pi: G\to \OO(\ell_p)$  a continuous
orthogonal  representation of $G$ on $\ell_p=\ell_p(X)$ for $1\leq p<\infty, p\neq2.$
   
By Banach's result, there exist mappings 
$$\vfi:G\to \CalS(X) \quad\text {and}\quad c: G\times X\to \S1$$
such that 
$$
\pi(g)(f)(x) = c(g^{-1},x) f(\vfi(g^{-1})(x)),
\tous   g\in G, f\in\ell_p(X),  x\in X.
$$
Since $\pi$ is a group homomorphism, one checks that $\vfi$
is also a group homomorphism; so, $\vfi$ defines an action of
$G$ on $X,$ which we hereafter denote simply by $(g, x)\mapsto gx.$
Moreover, $c: G\times X\to \S1$  satisfies the  cocycle relation $(*).$

Observe that   $\{\delta_x\ :\ x\in X\}$ is a discrete subset
of $\ell_p(X),$ equipped with the norm topology.  Since $\pi$ is continuous, it follows that  
the action of  $G$ on the discrete space $X$ is continuous.
Similarly,  one checks that $c: G\times X\to \S1$ is continuous.

In summary, to a continuous orthogonal representation $\pi$  
of $G$ on $\ell_p(X)$, $1<p<\infty, p\neq2,$ 
 is associated an  action of $G$ on  $X$
with open point stabilizers and a continuous cocycle  $c: G\times X\to \S1$.
(It is clear that, conversely, such an action of $G$ on $X$
and a continuous cocycle  $c: G\times X\to \S1$ define  
a continuous orthogonal representation of $G$ on $\ell_p(X).$) 

Set  $\pi^2=\la_{X}^c,$ where $\la_{X}^c$ 
is the permutation unitary representation of $G$ on $\ell_2(X)$
twisted by $c,$  as defined above.

Let  $\ell_p^{\pi(G)}$   and  $\ell_2^{\pi^2(G)}$ be 
the closed subspaces of $G$-invariant vectors in 
$\ell_p(X)$ and $\ell_2(X)$.
Let $\ell'_p(\pi)$ and  $\ell'_2(\pi^2)$ 
be the $G$-invariant complements of  $\ell_p^{\pi(G)}$   and  $\ell_2^{\pi^2(G)}$
as described in the introduction. (Observe that $\ell'_2(\pi^2)$ is the orthogonal
complement of  $\ell_2^{\pi^2(G)}$.)

The Mazur mappings 
 $M_{2,p}: \ell_2(X) \to \ell_p(X)$ and $M_{p,2}: \ell_p(X) \to \ell_2(X)$
 are the non-linear mappings defined  by 
 $$
 M_{2,p}(f)= (f/|f|) |f|^{2/p} \qquad\text{and} \qquad M_{p,2}(f)= (f/|f|) |f|^{p/2}.
 $$
 It is easily checked that $\pi^2(g)= M_{p,2} \pi(g) M_{2,p}$ and  $\pi(g)= M_{2,p} \pi^2(g) M_{p,2}$
 for all $g\in G.$
 The  mappings $M_{2,p}$ and $M_{2,p}$ are uniformly continuous  
 between the unit spheres of  $\ell_2$ and $\ell_p$ (\cite[Theorem~9.1]{BenLi}).
  % indeed,  the following inequalities hold $$  \frac{2}{p}\Vert f ? g\Vert_2\leq  \Vert M_{2,p}(f) ?M_{2,q}(y)kq ? C kx ? ykp/q
 As a consequence,  one obtains the   following crucial fact,
  established in Section 4.a  of \cite{BFGM}.
\begin{lemma}
\label{Lem1}
%\textbf{[BFGM]}
There  exists a sequence of almost  invariant 
vectors for the restriction of $\pi$ to $\ell'_p(\pi)$   
if and only if   there exists a  sequence of almost  invariant 
vectors for the restriction of $\pi^2$ to   $\ell'_2(\pi^2).$
$\bsq$ 
\end{lemma}

\section{Proof of Theorems~\ref{Theo1} and \ref{Theo2}}
 \label{S:Theo1-2}

\bigskip
\n
\textbf{Proof of Theorem~\ref{Theo1}}

\n
$(i) \Longrightarrow (ii):$ 

\n
Assume that $G$ does not have \Tlp.
 Then there exists an orthogonal representation  $\pi: G\to \OO( \ell_p)$   such that 
 the restriction of $\pi$ to $\ell'_p(\pi)$   has 
 a sequence of almost  invariant 
vectors.

 By Lemma~\ref{Lem1}, the 
unitary representation $\pi^2$ on $\ell_2$  associated to $\pi$
has  a sequence of almost  invariant 
vectors in $\ell'_2(\pi^2).$
On the other hand, Corollary~\ref{Cor-DirectSumMonoRep} shows that
$\pi^2 $ is unitary equivalent to a  direct sum 
of monomial representations associated to 
open subgroups. It follows that $1_G$
does not  have  Property $(T;  \R_{\rm mon}).$

\medskip
\n
$(ii) \Longrightarrow (i):$ 

\n
Assume that $1_G$
does not  have  Property $(T;  \R_{\rm mon}).$ 
Thus, there exists a family $(H_i, \chi_i)_{i\in I}$ of 
 open subgroups $H_i$ with  unitary characters $\chi_i$
 with the following property: 
 if $\pi^2$ denotes
 the representation  $\bigoplus_{i\in I}\Ind_{H_i}^{G} \chi_i$
 on $\H= \bigoplus_{i\in I}\ell_2(G/H_i),$
 there exists a sequence $(f_n)_n$ of almost invariant  
 vectors in   the orthogonal complement $\H'$ of the space
 $\H^{\pi(G)}$  of invariant vectors. 
 
For every $f\in \H,$ the projection of 
$f$ on $\ell_2(G/H_i)$  is non-zero for  at most countably
 many $i\in I.$  
 It follows that we can assume that the set $I$ is infinite countable.
 (If $I$ happens to be finite, we replace $I$ be $I\times \NN$ and set
 $H_{(i, n)}= H_i$.)
 
 Let  $ X=\coprod_{i\in I} G/H_i.$
  By Corollary~\ref{Cor-DirectSumMonoRep}, 
 $\pi^2$  is equivalent 
to the permutation representation $\la_{X}^c$  of $G$ on $\ell_2(X)$
twisted by a cocycle $c: G\times X\to \S1$. 
 
 We can associate to $\la_X^c$
the  orthogonal  representation $\pi: G\to \OO(\ell_p(X)),$
defined by the same formula.   Lemma~\ref{Lem1} 
shows that  $\pi$ has
 a sequence  of almost invariant  
vectors contained  in the $G$ -invariant complement $\ell'_p(\pi)$ of $\ell_p(X)^{\pi(G)}.$
Therefore, $G$ does not have \Tlp.
$\bsq$  

The  proof
of  Theorem~\ref{Theo2} will be an easy consequence of the following lemma.
\begin{lemma}
\label{Lem-MonoQuasiReg}
 Let $G$ be a totally disconnected group, $H$ an open subgroup
 and $\chi$ a continuous unitary character of $H.$ 
  There exists an open subgroup $L$ 
 of $G$ contained in $H$ such that the monomial representation
  $\Ind_{H}^{G} \chi$ is  weakly contained in 
  the quasi-regular representation $\la_{G/L}.$
  \end{lemma}
  \proof
  Since $G,$ and hence $H,$ is totally disconnected,
  every  neighbourhood of the group unit  in $H$ contains a compact open subgroup, 
   by  Dantzig's theorem (see, e.g., Theorem~7.7 in \cite{HewittRoss}).
  By continuity  of $\chi$, there exists a compact open subgroup
  $K$ of $H$ such that 
  $$|\chi(k)-1|<1 \tous k\in K.$$
   For every $k\in K,$ we then have 
  $|\chi(k)^n-1|<1$ for all  $n\in\NN$
 and hence $\chi(k)=1.$ Therefore $\chi$
 is trivial on $K.$
 
 Let $L$ be the   subgroup of $G$ generated
 by $K \cup [H, H].$
 Then $L$ is a normal and open subgroup of $H$ and $\chi$
 is trivial on $L$. So, $\chi$ factorizes to a 
 unitary character $\overline{\chi}$ of the abelian quotient group  $\overline{H}=H/L.$ 
 
 Since $\overline{H}$ is amenable,   $\overline{\chi}$ is weakly contained
 in the regular representation  $\la_{\overline{H}}$ of $\overline{H},$
 by the  Hulanicki-Reiter theorem (see Theorem~G.3.2 in \cite{BHV}).
 Hence, $\chi$ is weakly contained in the quasi-regular representation 
 $\la_{H/L},$ since  $\la_{H/L}=\la_{\overline{H}}\circ p,$
 where $p: H\to \overline{H}$ is the quotient homomorphism.
 By continuity of induction (see Theorem~F. 3.5 in \cite{BHV}),  it follows that
$\Ind_{H}^{G} \chi$ is  weakly contained in 
$$\Ind_{H}^G(\la_{H/L})\cong \la_{G/L}.\bsq$$

\bigskip
\n
\textbf{Proof of Theorem~\ref{Theo2}}

\n
In view of Theorem~\ref{Theo1}, it suffices to show that
Property $(T;  \R_{\rm quasi-reg}) $  implies 
$(T;  \R_{\rm mon}) $ for  second countable locally compact
and totally disconnected groups.

Assume that such a group $G$ does not have   $(T;  \R_{\rm mon}).$ 
Then there exists a family $(H_i, \chi_i)_{i\in I}$ of 
 open subgroups $H_i$ with  unitary characters $\chi_i$
 such that $1_G$ is weakly contained in the restriction
 $\pi'$ of 
 $$\pi:=\bigoplus_{i\in I}\Ind_{H_i}^{G} \chi_i$$
 to  the orthogonal complement
 of the $\pi(G)$-invariant vectors. 
 On the other hand, by Lemma~\ref{Lem-MonoQuasiReg}, there exists a family 
$(L_i)_{i\in I}$  of open subgroups $L_i$ of  $H_i$
such that $\pi$ is weakly 
contained in 
$$\rho:=\bigoplus_{i\in I}  \la_{G/L_i}.$$
This implies that $\pi'$ is weakly contained
in the restriction  of $\rho$
to the orthogonal complement
 of the $\rho(G)$-invariant vectors. 
 Hence,  $G$ does not have $(T;  \R_{\rm quasi-reg}).\bsq$

\section{Proof of  Theorem~\ref{Theo-HoweMoore}}
 \label{S:Consequences}

\n
(i) The proof is similar to Kazhdan's proof from \cite{Kazhdan} (see also Lemma 2.14 in \cite{Glasner-Monod}):
let $\C$ be the family of 
 open and compactly generated subgroups
 of  $G$. Since $G$ is locally compact, $1_G$ is weakly 
contained in the family of quasi-regular representations
$(\la_{G/H})_{H\in\C}.$ 
 Hence, by Theorem~\ref{Theo1},  
 there exists $H\in \C$ such that $G$ has a non-zero invariant vector
  in $\ell_2(G/H).$ 
 This implies that $H$ has finite index and therefore that $G$ is compactly generated.

\medskip
 \n
 (ii)  Assume, by contradiction, that  
  $G/\overline{[G,G]}$ is not compact.
  Then there exists a sequence $(\chi_n)_n$ of unitary characters
  of $G$ with $\chi_n\neq 1_G$ and such that $\lim_n\chi_n= 1_G$ uniformly on compact subsets of
  $G.$ This contradicts  Theorem~\ref{Theo1}.

 \medskip
 \n
 (iii) 
$\bullet$  Let $L$ be a finite index subgroup of $G$. We want to show that $L$ has \Tlp.
 
 Let $\L$  be the set of pairs $(H,\chi)$ consisting of an open subgroup
$H$ of $L$ and unitary character $\chi$ of $H.$ For $(H,\chi)\in \L,$
denote by $\la_{(H,\chi)}$
the induced representation $\Ind_H^L \chi.$ Set 
$$\rho=\bigoplus_{(H,\chi)\in \L}  \la_{(H,\chi)}.$$
Let  $\rho'$  be the restriction of $\rho$
 to the orthogonal complement of the space of $\rho(L)$-invariant vectors.

 Assume, by contradiction, that   $L$ does not have \Tlp\ for $p\neq 2.$
 Then, by Theorem~\ref{Theo1}, the trivial representation $1_L$ 
 of $L$ is  weakly contained in $\rho'.$
It follows, by continuity of induction, that 
 $ \la_{G/L}$ is weakly contained in $\Ind_L^G\rho',$ 
which is a subrepresentation of   
$$\bigoplus_{(H,\chi)\in\L}\Ind_L ^G (\la_{(H,\chi)})\cong \bigoplus_{(H,\chi)\in\L}\Ind_H^G\chi.$$

On the other hand,  $1_G$  contained in  $\la_{G/L},$ as $G/L$ is finite.
Therefore, $1_G$ is weakly contained in  $\Ind_L^G\rho'.$
However, $\Ind_L^G\rho'$ has no non-zero $G$-invariant vector, since 
$\rho'$ has no non-zero $L$-invariant ones (see \cite[Theorem E.3.1]{BHV}). 
This is a contradiction to Theorem~\ref{Theo1}. 
We conclude that $L$  has \Tlp\ for $p\neq 2.$

\bigskip
$\bullet$ Let $\widetilde{G}$ be a group containing $G$ as a subgroup of finite index.
We want to show that $\widetilde{G}$ has \Tlp.

Since   $G$  contains a normal subgroup of $\widetilde{G}$ 
of finite index and since this subgroup has  \Tlp\ for $p\neq 2,$
by what we have just seen above, we can assume  that $G$ is 
 a normal subgroup of $\widetilde{G}.$ 

Assume, by contradiction, that   $\widetilde{G}$ does not have \Tlp\ for $p\neq 2.$
Then there exists an orthogonal representation
 $\pi: \widetilde{G}\to \OO(\ell_p)$ which has a sequence  of almost  invariant 
vectors in the complement  $\ell'_p$  of $\pi(\widetilde{G})$-invariant vectors in $\ell_p.$

Let $\pi^2$ be the unitary representation of $\widetilde{G}$ on $\ell_2$ associated
to $\pi$, as in the beginning of this section. By Lemma~\ref{Lem1},
there exists a sequence $(\xi_n)_n$  of almost  invariant 
vectors in the orthogonal complement  $\ell'_2$  of $\pi^2(\widetilde{G})$-invariant vectors
in $\ell_2.$

Let $P:\ell_2' \to (\ell'_2)^{G}$ be the orthogonal projection on the
subspace of  $\pi^2(G)$-invariant vectors in $\ell'_2.$
Observe that $(\ell'_2)^G$ is invariant under $\pi^2(\widetilde{G}),$
since $G$ is normal in $\widetilde{G}.$

For every $n\in\NN,$ the vector  $\xi_n-P\xi_n$ belongs to  the orthogonal complement  of  $(\ell'_2)^G$
in $\ell'_2$. Hence, $\xi_n-P\xi_n$ belongs  to  the orthogonal complement  in $\ell_2$ of  
the space $(\ell_2)^G$ of $G$-invariant vectors, since $(\ell_2)^G= (\ell'_2)^G \oplus (\ell_2)^{\widetilde G}.$ 
Moreover, we have
$$\lim_n\Vert \pi^2(g)(\xi_n-P\xi_n)-(\xi_n-P\xi_n)\Vert= 0, \tous g\in G.
$$
It follows that $\inf_n \Vert \xi_n-P\xi_n\Vert=0;$ indeed, otherwise,
 $\frac{1}{\Vert\xi_n-P\xi_n\Vert}(\xi_n-P\xi_n)$ would be  
a sequence  of almost  invariant  
vectors in the orthogonal complement of $(\ell_2)^G$  in $\ell_2$ and, using Lemma~\ref{Lem1},
this would contradict the fact that $G$ has \Tlp.
Hence, upon passing to a subsequence, we can assume that   $\lim_n \Vert \xi_n-P\xi_n\Vert=0.$
%As a consequence, there exists $N$ such that $P\xi_n\neq 0$ for all $n\geq N.$

 Since $P\xi_n$  is $G$-invariant, we have can define the following sequence 
 $(\eta_n)_n$ of vectors in $\ell_2:$
$$
\eta_n= \dfrac{1}{\left|{\widetilde G}/G\right|}\sum_{t\in {\widetilde G}/G} \pi^2(t) P\xi_n.
$$
It is clear that  $\eta_n$ is $\pi^2(\widetilde G)$-invariant.  Moreover, we have
\begin{align*}
\Vert \eta_n-\xi_n\Vert &\leq \frac{1}{\left|{\widetilde G}/G\right|}\sum_{t\in {\widetilde G}/G} \Vert \pi^2(t) P\xi_n -\xi_n\Vert\\
&\leq  \frac{1}{\left|{\widetilde G}/G\right|}\sum_{t\in {\widetilde G}/G} \Vert \pi^2(t) P\xi_n -\pi^2(t)\xi_n\Vert+ \\
&\ + \frac{1}{\left|{\widetilde G}/G\right|}\sum_{t\in {\widetilde G}/G} \Vert \pi^2(t)\xi_n -\xi_n\Vert\\
&= \frac{1}{\left|{\widetilde G}/G\right|}\sum_{t\in {\widetilde G}/G} \Vert  P\xi_n -\xi_n\Vert
+ \frac{1}{\left|{\widetilde G}/G\right|}\sum_{t\in {\widetilde G}/G} \Vert \pi^2(t)\xi_n -\xi_n\Vert.
\end{align*}
It follows that $\lim_n\Vert \eta_n-\xi_n\Vert=0.$ Hence,  $\eta_n\neq 0 $ for sufficiently large $n,$
since $\Vert \xi_n\Vert=1.$

For every $t\in {\widetilde G},$ 
the vector  $\pi^2(t)( P\xi_n)$ belongs to $(\ell'_2)^G,$ 
since $(\ell'_2)^G$ is invariant under $\pi^2(\widetilde{G}).$
It follows that $\eta_n\in (\ell'_2)^G$ and, in particular,    $\eta_n\in \ell'_2.$ This is a contradiction,
as there are no non-zero $\pi^2(\widetilde G)$-invariant
vector in $\ell'_2.$

\medskip
 \n
 (iv) Since $G$ is totally disconnected,  we can find   a compact open subgroup $K$
 of $G, $ by van Dantzig's theorem. The amenability of  $G$ 
 implies  the amenability of its  action on $G/K:$
 $1_G$ is weakly contained in $\la_{G/K}$ (see Th\'eor\`eme on p. 28 in  \cite{Eymard}).
As $G$ has \Tlp, it  follows from Theorem~\ref{Theo2} that $G$
 has a non-zero invariant vector in $\ell_2(G/K).$ 
 Hence, $K$ has finite index in $G$ and  $G$ is compact.$\bsq$

\bigskip
%\section{Proof of  Theorem~\ref{Theo-HoweMoore}}
% \label{S:TotallyDisconnected}
\n
\textbf{Proof of  Theorem~\ref{Theo-HoweMoore}}
 
\n
Since $G$ has the Howe-Moore property, every proper open  subgroup of  $G$ 
 is compact (see \cite[Proposition 3.2]{CCLTV}).
 It follows that, for every proper open  subgroup $H$, the space $\ell_2(G/H)$ can be identified
with the  $G$-invariant subspace  of $L_2(G)$  of  functions on $G$ with are right $H$-invariant.
As a consequence, we see that  $\la_{G/H}$ is a subrepresentation of 
the regular representation $\la_G$.  
Denoting by $\L$  be the set of proper open subgroups of $G,$ 
this implies that   $\bigoplus_{H\in \L}  \la_{G/H}$ is  weakly contained in  the regular representation $\la_G$.
  
 On the other hand, since $G$ is not amenable, $1_G$ is not weakly contained
 in $\la_G,$ by the Reiter-Hulanicki theorem. It follows that 
 $1_G$ is not weakly contained
 in $\bigoplus_{H\in \L}  \la_{G/H}$ and
 Theorem~\ref{Theo2} shows that $G$ has \Tlp \ for $p\neq 2.\bsq$

 \section{Discrete groups with \Tlp}
 \label{S:Discrete}

The main tool for the proof  of Theorem~\ref{Theo3} is the following rigidity result
 concerning the unitary representations of 
 the groups appearing in the statement.
 
Given  locally compact second countable groups $G$ and $Q$ and a  continuous
homomorphism  $f : G\to Q$ with dense image, let us say as in \cite[Definition 4.2.2]{Cornulier}
that $f$ is a  \emph{resolution} if, whenever  a unitary representation $\pi$
 of $G$ has a sequence of almost invariant vectors, $\pi$ has
a non-zero subrepresentation $\rho$ 
 which factors through a unitary representation  $\widetilde \rho$ of $Q,$
 that is, $\rho= \widetilde{\rho}\circ f$.

  A weaker form  (which is however often sufficient for the applications) of 
  the following theorem  was established in \cite[Theorem 2.2]{Lubotzky-Zimmer};
  for  more general versions of it,
   see  \cite[Chap. III, Theorem 6.3]{Margulis}, \cite[Theorem A]{Bekka-Louvet}, and \cite[Theorem 4.3.1]{Cornulier}.
 \begin{theorem} 
 \label{Theo-LubotzkyZimmer}
 Let $G_1$ and $G_2$  be  second countable locally compact groups
 and $\Ga$  an irreducible lattice in 
 $G=G_1\times G_2.$ Assume that $G_1$ has \PT.
Then the canonical projection $p_2: \Ga \to G_2$ is a resolution.$\bsq$
\end{theorem}

%By a well-known result of B.H. Neumann (\cite[Lemma 4.1]{Neumann}),
%a group cannot be covered a finite union of cosets
%of subgroups of infinite index. 
The second ingredient for the proof of  Theorem~\ref{Theo3}
is the following elementary lemma.
% strengthening  of Neumann's theorem 
%in the case of groups which can be densely embedded in connected groups.

\begin{lemma}
\label{Lem-Tau}
Let $G$ be a connected topological group
and $\Ga$ a  dense subgroup of $G.$ 
Let $H_1, \dots, H_n$ be subgroups of $\Ga$
such that no $H_i$ is dense in $G.$
 Let $X= \bigcup_{i=1}^n a_i H_i$ be a  union
 of cosets of the $H_i$'s.
Then $\Ga\setminus X$
is dense in $G.$ 
\end{lemma}
\proof
 Assume, by contradiction, that ${\Ga\setminus X}$
 is not dense in $G.$
 Then there exists a non-empty open subset 
 $U$ of  $G$ such that  $U\cap (\Ga\setminus X)= \emptyset$
and hence
 $$U\cap (\overline{\Ga\setminus X})= \emptyset.$$
 Since $\Ga$ is dense in $G,$ it follows
 that $U$ is contained in 
 $$
 \overline{ X}= \bigcup_{i=1}^n a_i \overline H_i.
 $$
 This implies that some coset
 $a_i\overline{H_i}$ has non-empty interior in $G.$ Thus, 
$\overline{H_i}$ is open in $G$  and hence $\overline{H_i}=G,$
 since $G$ is connected. This is a contradiction. $\bsq$

 \bigskip
\n
 \textbf{Proof of  Theorem~\ref{Theo3}}

\n
Assume, by contradiction, that $\Ga$ does not have \Tlp. 
By  Theorem~\ref{Theo2}, there exists a family $\L$ of subgroups of $\Ga$
such that $1_\Ga$ is weakly contained
in the restriction $\pi'$ of  
$$\pi=\bigoplus_{H\in \L}  \la_{\Ga/H}$$ 
 to the orthogonal complement 
 the space of $\pi(\Ga)$-invariant vectors.

 Now,  $\Ga$ has \Ptau ; see Remark~\ref{Rem-Tau}. 
 So, we can assume that every $H\in \L$ has infinite index in $\Ga.$
 In particular, $\pi'$ coincides with $\pi$ and acts
 on $\H= \bigoplus_{H\in \L}  \ell_2(\Ga/H).$
 
   It follows from Theorem~\ref{Theo-LubotzkyZimmer} that
 there exists a  non-zero $\Ga$-invariant subspace
  $\K$ of $\H$
  such that the corresponding $\Ga$-representation factors through 
 a  (continuous) unitary representation  of $G_2.$
 
  Let $H\in\L$ be  such that 
 the orthogonal projection $P(\K)$ of $\K$ on $ \ell_2(\Ga/H)$ is non-zero.
 Then the restriction of $\pi$ to $P(\K)$  also factors through 
 a   unitary representation $\widetilde\pi$ of  $G_2$
 (see \cite[Lemma 4.1.5]{Cornulier}).

We claim that  $p_2(H)$ is not dense in $G_2.$
Indeed, assume that this is not the case.
The representation of $G,$ obtained by  inducing the restriction of $\pi$ 
to $P(\K),$ is  a subrepresentation of 
$$\Ind_\Ga^G(\la_{\Ga/H})\cong \la_{G/H}.$$
Since $\Ga$ is a lattice, it follows that there exists a non-zero $G_1$-invariant vector $f\in L_2(G/H).$
Lifting $f$ to $G,$ we obtain a measurable function $f: G\to \CCC$
such that, for every $g\in G_1,$  $f(g xh)= f(x)$ for almost every $x\in G$
and every $h\in H.$
Upon changing $f$ on a null set,  we can assume
that the previous equality holds for all $g\in G_1$, $x\in G$ and all $h\in H$
(see \cite[Lemma 2.2.16]{Zimmer}).
Thus,  $f$ is invariant under right translation by elements from $G_1H$. Since, by assumption, 
$G_1H$ is dense in $G,$
it follows that $f$ is constant, up to a null set (see \cite[Lemma 2.2.13]{Zimmer}). 
As $f\in L_2(G/H)$ is non zero, this implies that $G/H$  has finite volume. 
This is impossible, since $H$ has infinite index in $\Ga.$
So,   $p_2(H)$ is not dense in $G_2.$

 Let $f$ be a unit-vector in $P(\K).$ 
 We claim that we can find a sequence $\ga_n \in \Ga$ 
 with 
 $$\lim_n p_2(\ga_n) = e \quad\text{and}
 \quad \lim_n \langle \la_{\Ga/H}(\ga_n) f, f \rangle= 0.\eqno{(**)}$$
 This will yield the desired  contradiction, since  
 $$
\lim_n \langle \la_{\Ga/H}(\ga_n) f, f \rangle=\lim_n\langle \widetilde{\pi}(p_2(\ga_n)) f, f \rangle=1,
$$
by continuity of $\widetilde\pi.$

Let $\eps>0.$ There exists a finite subset $F$ of  $\Ga/H$ such that  
$$\sum_{x\notin F} |f(x)|^2 \leq  \eps.$$
Choose a set  $A\in  \Ga$  of representatives for the cosets $x\in F.$
Set
$$
X= \bigcup_{a,b\in A} p_2(b) p_2(H)p_2( a^{-1})= \bigcup_{a,b\in A} p_2(b a^{-1})p_2(H^{a}),
$$
which is a finite  union of cosets of subgroups conjugate to $p_2(H).$

As we have shown above, $p_2(H)$ is not dense in $G_2;$
the same is true for its conjugate subgroups.
Since $G_2$ is connected, it follows from Lemma~\ref{Lem-Tau}
that   $p_2(\Ga)\setminus X$ is dense in $G_2.$
We can therefore find a sequence $(\ga_n)_n$ in $\Ga$
with  $\lim_n p_2(\ga_n)= e$ such that, for all $n,$
 $$
 p_2(\ga_n) p_2(a) p_2(H) \cap p_2(b)p_2(H)=\emptyset, \tous a, b\in A,
 $$
 and hence
 $$
\ga_na H \cap b H=\emptyset, \tous a, b\in A, $$
for all $n.$
We then have 
\begin{align*}
\left|\langle \la_{\Ga/H}(\ga_n) f, f \rangle\right| &\leq \sum_{x\in F}  |f(\ga_n^{-1} x) f(x)| +
\sum_{x\notin F} |f(\ga_n^{-1}x) f(x)|\\
&\leq 2\Vert f\Vert  \left(\sum_{x\notin F} |f(x)|^2 \right)^{1/2}  \\
&\leq 2\sqrt{\eps}
\end{align*}
It follows that  the claim $(**)$ is satisfied by a subsequence of $(\ga_n)_n.\bsq$

\bigskip

\n
\textbf{Proof of  Proposition~\ref{Pro-F}}

\n
Assume that $\Ga$  does not have \Tlp.
Let $\L$ be the family of proper subgroups of  $\Ga$ and $X= \coprod_{H\in \L} \Ga/H.$ 
 Theorem~\ref{Theo2} shows that   the action of  $\Ga$ on  $X$ is amenable. 
 However, $\Ga$ has no fixed point in $X.$
Hence, $\Ga$ does not have  \PF. $\bsq$

\end{document}